\documentclass[12pt,thmsa]{article}%
\usepackage{amssymb}
\usepackage{amsfonts}
\usepackage{sw20bams}%
\usepackage{amsmath}%
\usepackage{graphicx}
\providecommand{\U}[1]{\protect\rule{.1in}{.1in}}
\begin{document}

\author{Steven R. Finch}
\title{Width Distributions for Convex Regular Polyhedra}
\date{March 12, 2016}
\maketitle

\begin{abstract}
The mean width is a measure on three-dimensional convex bodies that enjoys
equal status with volume and surface area \cite{Rt}. As the phrase suggests,
it is the mean of a probability density $f$. We verify formulas for mean
widths of the regular tetrahedron and the cube. Higher-order moments of
$f_{\text{tetra}}$ and $f_{\text{cube}}$ have not been examined until now.
\ Assume that each polyhedron has edges of unit length. We deduce that the
mean square width of the regular tetrahedron is $1/3+(3+\sqrt{3})/(3\pi)$ and
the mean square width of the cube is $1+4/\pi$.

\end{abstract}

\footnotetext{Copyright \copyright \ 2011, 2016 by Steven R. Finch. All rights
reserved.}Let $C$ be a convex body in $\mathbb{R}^{3}$. A\ \textbf{width} is
the distance between a pair of parallel $C$-supporting planes. Every unit
vector $u\in\mathbb{R}^{3}$ determines a unique such pair of planes orthogonal
to $u$ and hence a width $w(u)$. Let $u$ be uniformly distributed on the unit
sphere $S^{2}\subset$ $\mathbb{R}^{3}$. Then $w$ is a random variable and%
\[
\operatorname*{E}\left(  w_{\text{tetra}}\right)  =\frac{3}{2\pi}%
\arccos\left(  -\frac{1}{3}\right)
\]
for $C=$ the regular tetrahedron with edges of unit length and
\[
\operatorname*{E}\left(  w_{\text{cube}}\right)  =\frac{3}{2}
\]
for $C=$ the cube with edges of unit length \cite{Rt, F1, KR}. \ The
probability density of $w$ is not known. Our humble contribution is to verify
the preceding \textbf{mean width} results and to deduce the following
\textbf{mean square width} results:%
\[
\operatorname*{E}\left(  w_{\text{tetra}}^{2}\right)  =\frac{1}{3}\left(
1+\frac{3+\sqrt{3}}{\pi}\right)  ,
\]%
\[
\operatorname*{E}\left(  w_{\text{cube}}^{2}\right)  =1+\frac{4}{\pi}
\]
which appear to be new.

We start with two-dimensional analogs of these results in Sections 1 and 2.
\ The tetrahedral case (Section 3) is more difficult than the cubic case
(Section 4); details of the evaluation of a certain definite integral are
relegated to Section 5. \ We note that the phrases \textbf{mean breadth} and
\textbf{mean caliper diameter} are synonymous with mean width \cite{F1}, and
hope that this paper might inspire relevant computations for other convex bodies.

\section{Equilateral Triangle}

Given a convex region $C$ in $\mathbb{R}^{2}$, a\ width is the distance
between a pair of parallel $C$-supporting lines. Every unit vector
$u\in\mathbb{R}^{2}$ determines a unique such pair of lines orthogonal to $u$
and hence a width $w(u)$. Let $u$ be uniformly distributed on the unit circle
$S^{1}\subset$ $\mathbb{R}^{2}$. We wish to study the distribution of the
random variable $w$ in the case $C=$ the equilateral triangle with sides of
unit length.

For simplicity, let $\triangle$ be the equilateral triangle with vertices%
\[%
\begin{array}
[c]{ccccc}%
v_{1}=\left(  0,1\right)  , &  & v_{2}=\left(  \frac{\sqrt{3}}{2},-\frac{1}%
{2}\right)  , &  & v_{3}=\left(  -\frac{\sqrt{3}}{2},-\frac{1}{2}\right)  .
\end{array}
\]
At the end, it will be necessary to normalize by $\sqrt{3}$, the side-length
of $\triangle$.

Also let $\widetilde{\triangle}$ be the union of three overlapping disks of
radius $1/2$ centered at $v_{1}/2$, $v_{2}/2$, $v_{3}/2$. Clearly
$\triangle\subset\widetilde{\triangle}$ and $\widetilde{\triangle}$ has
centroid $(0,0)$. \ A \textbf{diameter} of $\widetilde{\triangle}$ is the
length of the intersection between $\widetilde{\triangle}$ and a line passing
through the origin.

Computing all widths of $\triangle$ is equivalent to computing all diameters
of $\widetilde{\triangle}$. The latter is achieved as follows. Fix a point
$(a,b)$ on the unit circle. \ The line $L$ passing through $(0,0) $ and
$(a,b)$ has parametric representation%
\[%
\begin{array}
[c]{ccccc}%
x=t\,a, &  & y=t\,b, &  & t\in\mathbb{R}%
\end{array}
\]
and hence $y=(b/a)x$ assuming $a\neq0$. \ The nontrivial intersection between
top circle and $L$ satisfies%
\[
x^{2}+\left(  \tfrac{b}{a}x-\tfrac{1}{2}\right)  ^{2}=\tfrac{1}{4}
\]
thus $x_{1}=a\,b$ since $a^{2}+b^{2}=1$. \ The nontrivial intersection between
right circle and $L$ satisfies%
\[
\left(  x-\tfrac{\sqrt{3}}{4}\right)  ^{2}+\left(  \tfrac{b}{a}x+\tfrac{1}%
{4}\right)  ^{2}=\tfrac{1}{4}
\]
thus $x_{2}=a\left(  \sqrt{3}a-b\right)  /2$ . The nontrivial intersection
between left circle and $L$ satisfies%
\[
\left(  x+\tfrac{\sqrt{3}}{4}\right)  ^{2}+\left(  \tfrac{b}{a}x+\tfrac{1}%
{4}\right)  ^{2}=\tfrac{1}{4}
\]
thus $x_{3}=-a\left(  \sqrt{3}a+b\right)  /2$ . \ 

We now examine all pairwise distances, squared, between the three intersection
points:%
\[
\left(  x_{i}-x_{j}\right)  ^{2}+\left(  \tfrac{b}{a}x_{i}-\tfrac{b}{a}%
x_{j}\right)  ^{2}=\left\{
\begin{array}
[c]{lll}%
\frac{1}{4}\left(  3-6\sqrt{3}a\,b+6b^{2}\right)  &  & \text{if }i=1,j=2\\
\frac{3}{4}\left(  1+2\sqrt{3}a\,b+2b^{2}\right)  &  & \text{if }i=1,j=3\\
3a^{2} &  & \text{if }i=2,j=3
\end{array}
\right.
\]
and define%
\[
g(a,b)=\max\left\{  \tfrac{1}{4}\left(  3-6\sqrt{3}a\,b+6b^{2}\right)
,\tfrac{3}{4}\left(  1+2\sqrt{3}a\,b+2b^{2}\right)  ,3a^{2}\right\}  .
\]
Therefore the mean width for $C$ is%
\[
\frac{1}{\sqrt{3}}\frac{1}{2\pi}%
{\displaystyle\int\limits_{0}^{2\pi}}
\sqrt{g(\cos\theta,\sin\theta)}\,d\theta=\frac{3}{\pi}=\frac{\text{perimeter
of }C}{\pi}
\]
and the mean square width is%
\[
\frac{1}{3}\frac{1}{2\pi}%
{\displaystyle\int\limits_{0}^{2\pi}}
g(\cos\theta,\sin\theta)\,d\theta=\frac{1}{2}\left(  1+\frac{3\sqrt{3}}{2\pi
}\right)  .
\]

\section{Square}

We turn to the case $C=$ the square with sides of unit length. For simplicity,
let $\square$ be the square with vertices%
\[%
\begin{array}
[c]{ccccccc}%
v_{1}=\left(  \frac{\sqrt{2}}{2},\frac{\sqrt{2}}{2}\right)  , &  &
v_{2}=\left(  \frac{\sqrt{2}}{2},-\frac{\sqrt{2}}{2}\right)  , &  &
v_{3}=\left(  -\frac{\sqrt{2}}{2},\frac{\sqrt{2}}{2}\right)  , &  &
v_{4}=\left(  -\frac{\sqrt{2}}{2},-\frac{\sqrt{2}}{2}\right)  .
\end{array}
\]
At the end, it will be necessary to normalize by $\sqrt{2}$, the side-length
of $\square$.

Also let $\widetilde{\square}$ be the union of four overlapping disks of
radius $1/2$ centered at $v_{1}/2$, $v_{2}/2$, $v_{3}/2$, $v_{4}/2$. A
diameter of $\widetilde{\square}$ is the length of the intersection between
$\widetilde{\square}$ and a line passing through the origin.

Computing all widths of $\square$ is equivalent to computing all diameters of
$\widetilde{\square}$. Fix a point $(a,b)$ on the unit circle. \ The line $L$
passing through $(0,0)$ and $(a,b)$ can be represented as $y=(b/a)x$ assuming
$a\neq0$. \ The nontrivial intersection between northeast circle and $L$
satisfies%
\[
\left(  x-\tfrac{\sqrt{2}}{4}\right)  ^{2}+\left(  \tfrac{b}{a}x-\tfrac
{\sqrt{2}}{4}\right)  ^{2}=\tfrac{1}{4}
\]
thus $x_{1}=\sqrt{2}a(a+b)/2$ since $a^{2}+b^{2}=1$. \ The nontrivial
intersection between southeast circle and $L$ satisfies%
\[
\left(  x-\tfrac{\sqrt{2}}{4}\right)  ^{2}+\left(  \tfrac{b}{a}x+\tfrac
{\sqrt{2}}{4}\right)  ^{2}=\tfrac{1}{4}
\]
thus $x_{2}=\sqrt{2}a(a-b)/2$. The nontrivial intersection between northwest
circle and $L$ satisfies%
\[
\left(  x+\tfrac{\sqrt{2}}{4}\right)  ^{2}+\left(  \tfrac{b}{a}x-\tfrac
{\sqrt{2}}{4}\right)  ^{2}=\tfrac{1}{4}
\]
thus $x_{3}=-\sqrt{2}a(a-b)/2$. The nontrivial intersection between southwest
circle and $L$ satisfies%
\[
\left(  x+\tfrac{\sqrt{2}}{4}\right)  ^{2}+\left(  \tfrac{b}{a}x+\tfrac
{\sqrt{2}}{4}\right)  ^{2}=\tfrac{1}{4}
\]
thus $x_{4}=-\sqrt{2}a(a+b)/2$.

We now examine all pairwise distances, squared, between the four intersection
points:%
\[
\left(  x_{i}-x_{j}\right)  ^{2}+\left(  \tfrac{b}{a}x_{i}-\tfrac{b}{a}%
x_{j}\right)  ^{2}=\left\{
\begin{array}
[c]{lll}%
2b^{2} &  & \text{if }i=1,j=2\\
2a^{2} &  & \text{if }i=1,j=3\\
2+4a\,b &  & \text{if }i=1,j=4\\
2-4a\,b &  & \text{if }i=2,j=3\\
2a^{2} &  & \text{if }i=2,j=4\\
2b^{2} &  & \text{if }i=3,j=4
\end{array}
\right.
\]
and define%
\begin{align*}
g(a,b)  & =\max\left\{  2a^{2},2b^{2},2(1+2a\,b),2(1-2a\,b)\right\} \\
& =\max\left\{  2(1+2a\,b),2(1-2a\,b)\right\}  .
\end{align*}
Therefore the mean width for $C$ is%
\[
\frac{1}{\sqrt{2}}\frac{1}{2\pi}%
{\displaystyle\int\limits_{0}^{2\pi}}
\sqrt{g(\cos\theta,\sin\theta)}\,d\theta=\frac{4}{\pi}=\frac{\text{perimeter
of }C}{\pi}
\]
and the mean square width is%
\[
\frac{1}{2}\frac{1}{2\pi}%
{\displaystyle\int\limits_{0}^{2\pi}}
g(\cos\theta,\sin\theta)\,d\theta=1+\frac{2}{\pi}.
\]

\section{Regular Tetrahedron}

Returning to the three-dimensional setting of the introduction, let $C=$ the
regular tetrahedron with edges of unit length.

For simplicity, let $\triangle$ be the tetrahedron with vertices%
\[%
\begin{array}
[c]{lll}%
v_{1}=\left(  0,0,1\right)  , &  & v_{2}=\left(  \frac{2\sqrt{2}}{3}%
,0,-\frac{1}{3}\right)  ,\\
v_{3}=\left(  -\frac{\sqrt{2}}{3},\sqrt{\frac{2}{3}},-\frac{1}{3}\right)  , &
& v_{4}=\left(  -\frac{\sqrt{2}}{3},-\sqrt{\frac{2}{3}},-\frac{1}{3}\right)  ,
\end{array}
\]
At the end, it will be necessary to normalize by $2\sqrt{2/3}$, the
edge-length of $\triangle$.

Also let $\widetilde{\triangle}$ be the union of four overlapping balls of
radius $1/2$ centered at $v_{1}/2$, $v_{2}/2$, $v_{3}/2$, $v_{4}/2$. Clearly
$\triangle\subset\widetilde{\triangle}$ and $\widetilde{\triangle}$ has
centroid $(0,0,0)$. \ A diameter of $\widetilde{\triangle}$ is the length of
the intersection between $\widetilde{\triangle}$ and a line passing through
the origin.

Computing all widths of $\triangle$ is equivalent to computing all diameters
of $\widetilde{\triangle}$. The latter is achieved as follows. Fix a point
$(a,b,c)$ on the unit sphere. \ The line $L$ passing through $(0,0,0)$ and
$(a,b,c)$ has parametric representation%
\[%
\begin{array}
[c]{ccccccc}%
x=t\,a, &  & y=t\,b, &  & z=t\,c, &  & t\in\mathbb{R}%
\end{array}
\]
and hence $y=(b/a)x$, $z=(c/a)x$ assuming $a\neq0$. \ The nontrivial
intersection between top sphere and $L$ satisfies%
\[
x^{2}+\left(  \tfrac{b}{a}x\right)  ^{2}+\left(  \tfrac{c}{a}x-\tfrac{1}%
{2}\right)  ^{2}=\tfrac{1}{4}
\]
thus $x_{1}=a\,c$ since $a^{2}+b^{2}+c^{2}=1$. \ The nontrivial intersection
between front sphere and $L$ satisfies%
\[
\left(  x-\tfrac{\sqrt{2}}{3}\right)  ^{2}+\left(  \tfrac{b}{a}x\right)
^{2}+\left(  \tfrac{c}{a}x+\tfrac{1}{6}\right)  ^{2}=\tfrac{1}{4}
\]
thus $x_{2}=a\left(  2\sqrt{2}a-c\right)  /3$ . \ The nontrivial intersection
between left sphere and $L$ satisfies%

\[
\left(  x+\tfrac{\sqrt{2}}{6}\right)  ^{2}+\left(  \tfrac{b}{a}x+\tfrac
{\sqrt{6}}{6}\right)  ^{2}+\left(  \tfrac{c}{a}x+\tfrac{1}{6}\right)
^{2}=\tfrac{1}{4}
\]
thus $x_{3}=-a\left(  \sqrt{2}a+\sqrt{6}b+c\right)  /3$. \ The nontrivial
intersection between right sphere and $L$ satisfies%

\[
\left(  x+\tfrac{\sqrt{2}}{6}\right)  ^{2}+\left(  \tfrac{b}{a}x-\tfrac
{\sqrt{6}}{6}\right)  ^{2}+\left(  \tfrac{c}{a}x+\tfrac{1}{6}\right)
^{2}=\tfrac{1}{4}
\]
thus $x_{4}=a\left(  -\sqrt{2}a+\sqrt{6}b-c\right)  /3$. \ 

We now examine all pairwise distances, squared, between the four intersection
points:%
\begin{align*}
& \left(  x_{i}-x_{j}\right)  ^{2}+\left(  \tfrac{b}{a}x_{i}-\tfrac{b}{a}%
x_{j}\right)  ^{2}+\left(  \tfrac{c}{a}x_{i}-\tfrac{c}{a}x_{j}\right)  ^{2}\\
& =\left\{
\begin{array}
[c]{lll}%
\frac{8}{9}\left(  a^{2}-2\sqrt{2}a\,c+2c^{2}\right)  &  & \text{if }i=1,j=2\\
\tfrac{1}{9}\left(  \sqrt{2}a+\sqrt{6}b+4c\right)  ^{2} &  & \text{if
}i=1,j=3\\
\tfrac{1}{9}\left(  \sqrt{2}a-\sqrt{6}b+4c\right)  ^{2} &  & \text{if
}i=1,j=4\\
\frac{2}{3}\left(  3a^{2}+2\sqrt{3}a\,b+b^{2}\right)  &  & \text{if }i=2,j=3\\
\frac{2}{3}\left(  3a^{2}-2\sqrt{3}a\,b+b^{2}\right)  &  & \text{if }i=2,j=4\\
\frac{8}{3}b^{2} &  & \text{if }i=3,j=4
\end{array}
\right.
\end{align*}
and define%
\begin{align*}
g(a,b)  & =\max\left\{  \tfrac{8}{9}\left(  a^{2}-2\sqrt{2}a\,c+2c^{2}\right)
,\tfrac{1}{9}\left(  \sqrt{2}a+\sqrt{6}b+4c\right)  ^{2},\tfrac{1}{9}\left(
\sqrt{2}a-\sqrt{6}b+4c\right)  ^{2},\right. \\
& \left.  \;\;\;\;\;\;\;\;\;\tfrac{2}{3}\left(  3a^{2}+2\sqrt{3}%
a\,b+b^{2}\right)  ,\tfrac{2}{3}\left(  3a^{2}-2\sqrt{3}a\,b+b^{2}\right)
,\tfrac{8}{3}b^{2}\right\}  .
\end{align*}
As for the equilateral triangle, no simplication of $g$ seems possible. \ The
mean width for $C$ is%
\[
\frac{1}{2\sqrt{2/3}}\frac{1}{4\pi}%
{\displaystyle\int\limits_{0}^{2\pi}}
{\displaystyle\int\limits_{0}^{\pi}}
\sqrt{g(\cos\theta\sin\varphi,\sin\theta\sin\varphi,\cos\varphi)}\sin
\varphi\,d\varphi\,d\theta=\frac{3}{2\pi}\arccos\left(  -\frac{1}{3}\right)
\]
and the mean square width is%
\[
\frac{1}{8/3}\frac{1}{4\pi}%
{\displaystyle\int\limits_{0}^{2\pi}}
{\displaystyle\int\limits_{0}^{\pi}}
g(\cos\theta\sin\varphi,\sin\theta\sin\varphi,\cos\varphi)\sin\varphi
\,d\varphi\,d\theta=\frac{1}{3}\left(  1+\frac{3+\sqrt{3}}{\pi}\right)  .
\]
Details on the final integral are given in Section 5.

\section{Cube}

We turn to the case $C=$ the cube with edges of unit length. For simplicity,
let $\square$ be the cube with vertices%
\[
v_{k}=\left(  \pm\tfrac{\sqrt{3}}{3},\pm\tfrac{\sqrt{3}}{3},\pm\tfrac{\sqrt
{3}}{3}\right)
\]
for $1\leq k\leq8$. \ At the end, it will be necessary to normalize by
$2/\sqrt{3}$, the edge-length of $\square$.

Also let $\widetilde{\square}$ be the union of eight overlapping balls of
radius $1/2$ centered at $v_{1}/2$, $v_{2}/2$, $v_{3}/2$, \ldots, $v_{8}/2$. A
diameter of $\widetilde{\square}$ is the length of the intersection between
$\widetilde{\square}$ and a line passing through the origin.

Computing all widths of $\square$ is equivalent to computing all diameters of
$\widetilde{\square}$. Fix a point $(a,b,c)$ on the unit sphere. \ The line
$L$ passing through $(0,0,0)$ and $(a,b,c)$ can be represented as $y=(b/a)x$,
$z=(c/a)x$ assuming $a\neq0$. \ The nontrivial intersection between each of
the eight spheres and $L$ satisfies%
\[
\left(  x\pm\tfrac{\sqrt{3}}{6}\right)  ^{2}+\left(  \tfrac{b}{a}x\pm
\tfrac{\sqrt{3}}{6}\right)  ^{2}+\left(  \tfrac{c}{a}x\pm\tfrac{\sqrt{3}}%
{6}\right)  ^{2}=\tfrac{1}{4}
\]
and consequently%
\begin{align*}
g(a,b)  & =\max\left\{  \tfrac{4}{3}a^{2},\tfrac{4}{3}b^{2},\tfrac{4}{3}%
c^{2},\tfrac{4}{3}\left(  1+2a\,b-c^{2}\right)  ,\tfrac{4}{3}\left(
1+2a\,c-b^{2}\right)  ,\right. \\
& \;\;\;\;\;\;\;\;\;\tfrac{4}{3}\left(  1-2a\,b-c^{2}\right)  ,\tfrac{4}%
{3}\left(  1-2a\,c-b^{2}\right)  ,\tfrac{4}{3}\left(  b+c\right)  ^{2}%
,\tfrac{4}{3}\left(  b-c\right)  ^{2},\\
& \;\;\;\;\;\;\;\;\;\tfrac{4}{3}\left(  1+2a\,b+2a\,c+2b\,c\right)  ,\tfrac
{4}{3}\left(  1+2a\,b-2a\,c-2b\,c\right)  ,\\
& \left.  \;\;\;\;\;\;\;\;\,\tfrac{4}{3}\left(  1-2a\,b-2a\,c+2b\,c\right)
,\tfrac{4}{3}\left(  1-2a\,b+2a\,c-2b\,c\right)  \right\}
\end{align*}
after examining all 28 pairwise distances and extracting 13 distinct
expressions. As for the square (in which $g$ simplified to a maximum over two
terms), here $g$ simplifies to a maximum over four terms:
\begin{align*}
g(a,b)  & =\max\left\{  \tfrac{4}{3}\left(  1+2a\,b+2a\,c+2b\,c\right)
,\tfrac{4}{3}\left(  1+2a\,b-2a\,c-2b\,c\right)  ,\right. \\
& \left.  \;\;\;\;\;\;\;\;\;\tfrac{4}{3}\left(  1-2a\,b-2a\,c+2b\,c\right)
,\tfrac{4}{3}\left(  1-2a\,b+2a\,c-2b\,c\right)  \right\}  .
\end{align*}
The mean width for $C$ is%
\[
\frac{1}{2/\sqrt{3}}\frac{1}{4\pi}%
{\displaystyle\int\limits_{0}^{2\pi}}
{\displaystyle\int\limits_{0}^{\pi}}
\sqrt{g(\cos\theta\sin\varphi,\sin\theta\sin\varphi,\cos\varphi)}\sin
\varphi\,d\varphi\,d\theta=\frac{3}{2}
\]
and the mean square width is%
\[
\frac{1}{4/3}\frac{1}{4\pi}%
{\displaystyle\int\limits_{0}^{2\pi}}
{\displaystyle\int\limits_{0}^{\pi}}
g(\cos\theta\sin\varphi,\sin\theta\sin\varphi,\cos\varphi)\sin\varphi
\,d\varphi\,d\theta=1+\frac{4}{\pi}.
\]

\section{A\ Definite Integral}

Considerable work is required to evaluate the definite integral at the end of
Section 3. A plot of the surface
\[
(\theta,\varphi)\longmapsto\sqrt{\frac{g(\cos\theta\sin\varphi,\sin\theta
\sin\varphi,\cos\varphi)}{8/3}}
\]
appears in Figure 1, where $0\leq\theta\leq2\pi$ and $0\leq\varphi\leq\pi$.
\ Figure 2 contains the same surface, but viewed from above. \ Our focus will
be on the part of the surface in the lower right corner, specifically
$0\leq\theta\leq\pi/3$ and $2\lessapprox\varphi\leq\pi$. \ The volume under
this part is $1/24^{\text{th}}$ of the volume under the full surface.

We need to find the precise lower bound on $\varphi$ as a function of $\theta
$. Recall the formula for\ $g$ as a maximum over six terms in Section 3; let
$g_{\ell}$ denote the $\ell^{\text{th}}$ term, where $1\leq\ell\leq6$. \ Then
the lower bound on $\varphi$ is found by solving the equation%
\[
g_{1}(\cos\theta\sin\varphi,\sin\theta\sin\varphi,\cos\varphi)=g_{4}%
(\cos\theta\sin\varphi,\sin\theta\sin\varphi,\cos\varphi)
\]
for $\varphi$. We obtain $\varphi(\theta)=2\arctan(h(\theta))$, where
\[
h(\theta)=\frac{\cos\theta+\sqrt{3}\sin\theta+\sqrt{10-\cos(2\theta)+\sqrt
{3}\sin(2\theta)}}{2\sqrt{2}}
\]
and, in particular,
\[
\varphi(0)=2\arctan\left(  \sqrt{2}\right)  \approx1.9106,
\]%
\[
\varphi(\pi/3)=2\arctan\left(  \left(  1+\sqrt{3}\right)  /\sqrt{2}\right)
\approx2.1862.
\]
It follows that $g=g_{1}$ for $0\leq\theta\leq\pi/3$ and $2\arctan
(h)\leq\varphi\leq\pi$. \ Now we have%
\begin{align*}
& \frac{1}{8/3}\frac{1}{4\pi}%
{\displaystyle\int}
g_{1}(\cos\theta\sin\varphi,\sin\theta\sin\varphi,\cos\varphi)\sin
\varphi\,d\varphi\\
& =\frac{\left(  -3+\cos(2\theta)\right)  \cos(3\varphi)-3\left(
7+3\cos(2\theta)\right)  \cos\varphi-16\sqrt{2}\cos\theta\sin^{3}\varphi
}{288\pi}%
\end{align*}
and
\[
\left.  \cos(3\varphi)\right\vert _{2\arctan(h)}^{\pi}=-2\frac{\left(
1-3h^{2}\right)  ^{2}}{(1+h^{2})^{3}},
\]%
\[%
\begin{array}
[c]{ccc}%
\left.  \cos(\varphi)\right\vert _{2\arctan(h)}^{\pi}=-\dfrac{2}{1+h^{2}}, &
& \left.  \sin(\varphi)\right\vert _{2\arctan(h)}^{\pi}=-\dfrac{2h}{1+h^{2}}%
\end{array}
\]
therefore
\begin{align*}
& \frac{1}{8/3}\frac{1}{4\pi}%
{\displaystyle\int_{2\arctan(h)}^{\pi}}
g_{1}(\cos\theta\sin\varphi,\sin\theta\sin\varphi,\cos\varphi)\sin
\varphi\,d\varphi\\
& =\frac{6h^{4}+8\sqrt{2}h^{3}\cos(\theta)+3h^{2}\left(  1+\cos(2\theta
)\right)  +\left(  3+\cos(2\theta)\right)  }{18\pi(1+h^{2})^{3}}.
\end{align*}
Integrating this expression from $0$ to $\pi/3$ gives the desired formula for
$\operatorname*{E}\left(  w_{\text{tetra}}^{2}\right)  $.

\bigskip

\bigskip%
\begin{figure}[ptb]%
\centering
\includegraphics[
height=4.1131in,
width=5.2217in
]%
{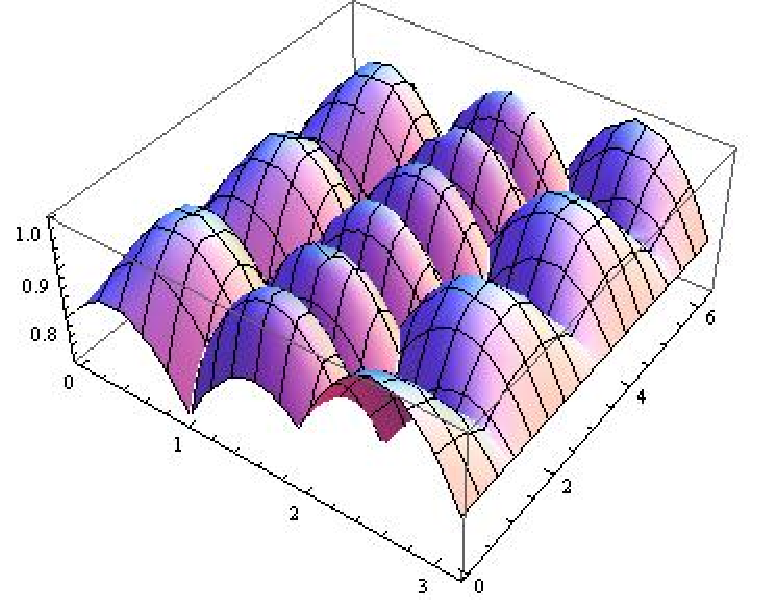}%
\caption{Surface plot of $\sqrt{3g/8}$, where $0\leq\theta\leq2\pi$ and
$0\leq\varphi\leq\pi$.}%
\end{figure}

\bigskip%

\begin{figure}[ptb]%
\centering
\includegraphics[
height=4.1131in,
width=5.2217in
]%
{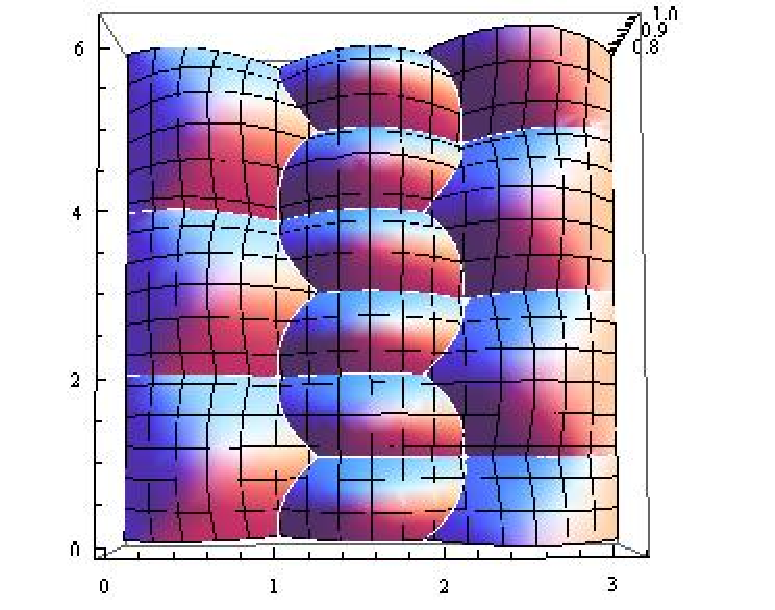}%
\caption{Another view of $\sqrt{3g/8}$, with contours of intersection.}%
\end{figure}

\end{document}